\documentclass[12pt]{article}
\usepackage[standard]{arvin}
\def\C{\mathcal{C}}
\DeclareMathOperator{\wt}{wt}

\DeclareMathOperator{\amax}{argmax}
\def\eps{\epsilon}
\newcommand{\ent}{H}
\usepackage{authblk}
\title{Unweighted One-Sided Code Sparsifiers and Thin Subgraphs}
\author[1]{Shayan Oveis Gharan}
 \affil{\small University of Washington, \textsf{shayan@cs.washington.edu}}
\author[2]{Arvin Sahami}
\affil{\small University of British Columbia, \textsf{arvin52@student.ubc.ca}}
\date{February 4, 2025}

\AddToHookNext{shipout/background}{
 \begin{tikzpicture}[remember picture, overlay,inner sep=0pt,outer sep=0pt,opacity=1]
 \node[anchor=north east] at ([xshift=2.7cm,yshift=3cm]current page.north east) {
  \includegraphics[scale = 0.0175]{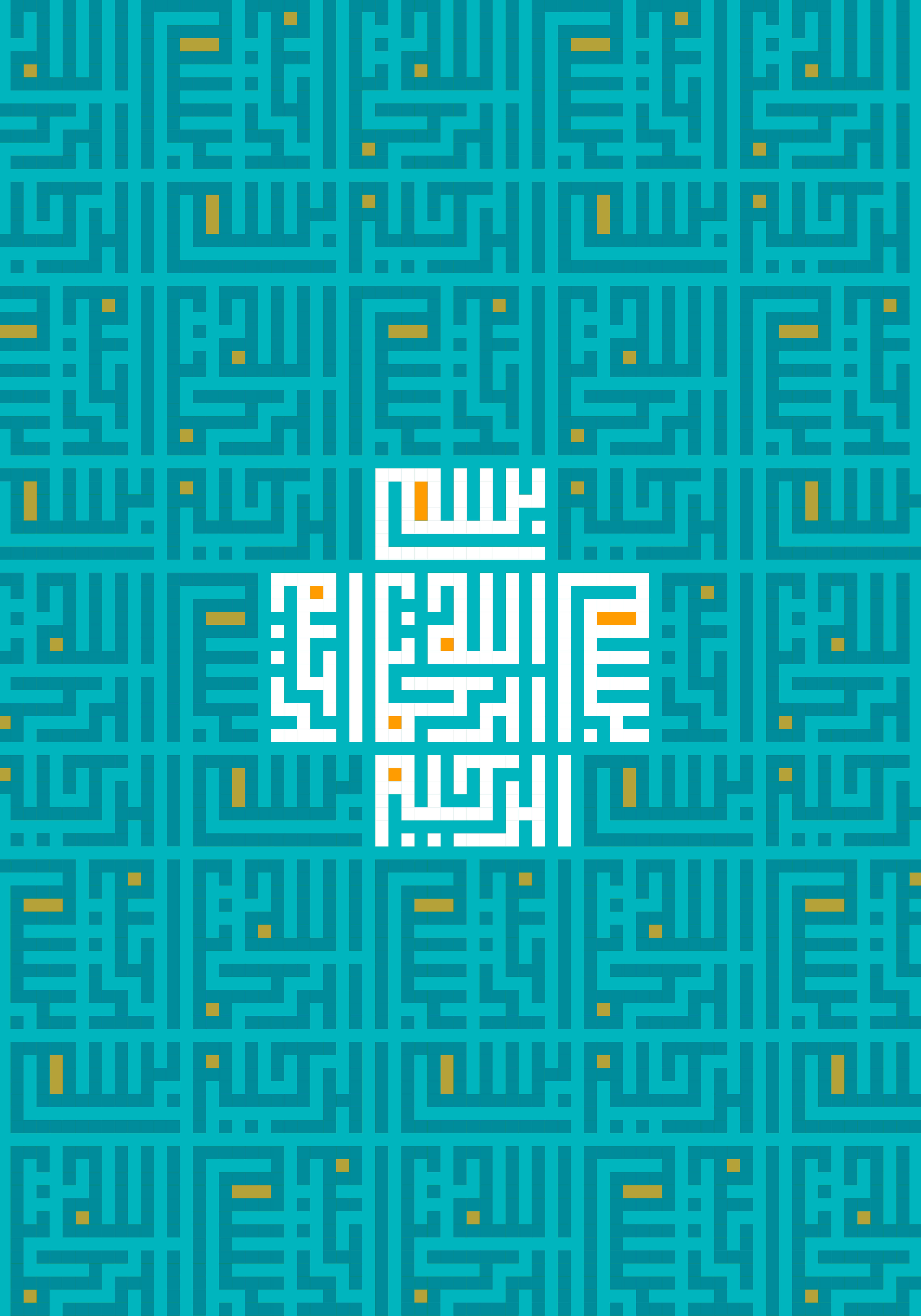}
 };
 \end{tikzpicture}
}

\newcommand{\cd}{k} %code dimension
\newcommand{\gk}{d} %changing graph connectivity symbol from k to d 
\newcommand{\size}{\abs}

\begin{document}
\maketitle
\begin{abstract}
    For a linear code $\cC \subseteq \bF_2^n$ and $\alpha \in [0,1]$, call a set $S \subseteq [n]$ an (unweighted) one-sided $\alpha$-sparsifier of $\cC$ if for all $c \in \cC$, $\wt(c_S)\geq \alpha \cdot \wt(c)$, where $c_S$ is the projection of $c$ onto the coordinates in $S$ and $\wt(c)$ is the Hamming weight of $c$.
    \\
    We show that every $k$-dimensional linear code $\cC\subseteq \bF_2^n$ has at least $2^{n - k}$ many unweighted one-sided  $1/2$-sparsifiers and hence one of size at most $n/2 + O(\sqrt{n k})$.
    As an application, letting $\cC \subseteq \bF_2^E$ denote the cut-space of a graph $G=\pair V E$, we show a lower bound of $2^{\size E - (\size V - 1)}$ on the number of $1/2$-thin subgraphs of $G$ and the existence of a $1/2$-thin subgraph with at least $\size E/2-O(\sqrt{\size E \cdot \size V})$ edges.
    \\
    In contrast to previous results on thin subgraphs, our proofs are purely "combinatorial".
\end{abstract}
\section{Introduction}
Recently, Khanna-Putterman-Sudan \cite{KPS24}  initiated the study of code sparsification where they proved that any linear code $\cC \subset \bF_2^n$ of dimension $\cd$ has a ``{\bf weighted}'' two-sided code sparsifier of size $\widetilde{O}(\cd/\eps^2)$. That is, they show there exists a set $S\subseteq [n]$ with $\size S \leq \widetilde{O}(\cd/\eps^2)$ together with non-negative weights $w \colon S\to\R_{\geq 0}$,  such that  for any $c\in \cC$, 
\[
    (1-\eps)\wt(c)\leq \sum_{i\in S} w_i c_i \leq(1+\eps)\wt(c)
\]
where $\wt(c)$ denotes the hamming weight of $c$ (see also  \cite{BG25} for extensions).
In the rest of the paper, for a codeword $c \in \bF_2^n$ and a set $S \subseteq [n]$, $c_s$ refers to the projection of $c$ onto the coordinates appearing in $S$.
\\
In our main theorem, we show that every linear code $\cC\subseteq \mathbb{F}_2^n$ has an {\em unweighted} one-sided sparsifier $S$ such that for any $c\in \C$, $\wt(c_S)\geq \wt(c)/2$ where $\abs{S}\approx n/2$. In other words, (assuming $\cd\ll n$), we can increase the rate of the code by a factor of $2$ while preserving almost the same relative distance.

We begin by defining our main object of study. 
\begin{defn}
    Let $\cC \subseteq \bF_2^n$ be a linear code (that is, $x,y\in\C$ implies $x+y\in \C$).
    For $\alpha \in [0, 1]$  call a set $S \subseteq [n]$ an unweighted one-sided $\alpha$-sparsifier of $\cC$ if for all $c \in \cC$ we have 
    \[
        \wt(c_S) \geq \alpha \wt(c).
    \]
\end{defn}
This definition provides an abstraction for the notion of thinness in graph theory which we discuss later. 

Our main theorem is the following.
\begin{thm}[Main]\label{thm:main}
    Any $\cd$-dimensional linear code $\cC \subseteq \bF_2^n$ has at least $2^{n-k}$ different unweighted one-sided $1/2$-sparsifiers. In other words, a uniformly random $S\subseteq [n]$ is an unweighted $1/2$-sparsifier w.p. at least $2^{-k}$.
\end{thm}
Let $\ent \colon [0, 1]\to \R$ be the entropy function defined as
\[
    \ent(x)=
    x\log_2 \frac1x + (1-x)\log_2\frac1{1-x}.
\]
\Cref{thm:main} immediately yields the following corollary, which we think is worth stating separately.
\begin{cor}\label{cor:smallS}
    Any $\cd$-dimensional linear code $\cC \subseteq \bF_2^n$ as above, has an unweighted one-sided $1/2$-sparsifier $S \subseteq [n]$ with $\size S \leq n(1/2 + \eps)$ where $\eps \in \open 0 {1/2}$ satisfies $\ent(1/2 - \eps) < 1-k/n$. In particular,  $\eps = \sqrt{\frac{\ln 2}{2}\cdot\frac{k}{n}}$ works.
\end{cor}
\begin{proof}[proof of \cref{cor:smallS}]
    We recall that for any $\gamma \in [0,1/2]$ we have
    \[
        \sum_{j=0}^{\gamma n} {n\choose j} \leq 2^{n\cdot \ent(\gamma)}.
    \]
    Note that the number of sets in $2^{[n]}$ with at least $n(1/2 + \eps)$ elements is 
    \[
        \sum_{j=0}^{(1/2 - \eps) n} {n\choose j} \leq 2^{n\cdot \ent(1/2 - \eps)}
        < 2 ^ {n - k}.
    \]
    But by \cref{thm:main} we know there are at least $2^{n-k}$ unweighted one-sided $1/2$-sparsifiers in $2^{[n]}$, thus at least one of them has at most $n(1/2 + \eps)$ elements. 
    \\
    The fact that $\eps = \sqrt{\frac{\ln 2}{2}\frac{k}{n}}$ works follows from the inequality $\ent(1/2 - x) \leq 1 - \frac{2}{\ln 2} x^2$ (with equality only at $x = 0$) which itself follows from the Taylor expansion of $H(x)$ at $1/2$. 
\end{proof}

\newcommand{\const}{\sqrt{\ln 2} (1 + \sqrt 2)}
We note that one can construct unweighted one-sided $\alpha$-sparsifiers for $\alpha > 1/2$ by a repeated application of \cref{cor:smallS}.
\begin{cor} \label{cor:bigAlpha}
    Let $\cC$ be as in \cref{cor:smallS}. Let $\alpha = 1 - \frac{1}{2^\ell}$ for an arbitrary positive integer $\ell$. Then $\cC$ has an unweighted one-sided $\alpha$-sparsifier $S$ with $\size S \leq \alpha n + c \sqrt{nk}$ for $c= \const$.
\end{cor}
We leave the proof to the appendix.
\\
We conclude this part with a few remarks.
\begin{itemize} 
\item The above bound is almost tight for any non-degenerate code $\C$ (i.e., one that for any coordinate $i$ has a codeword $c$ with $c_i\neq 0$). In particular any unweighted one-sided $\alpha$-sparsifier of such a code has at least $\alpha n$ coordinates. 
\item Our result is incomparable to that of \cite{KPS24}. This is because although we only guarantee a one-sided sparsifier, our set has a much simpler description since it is unweighted.
Furthermore, unlike \cite{KPS24} the size of the sparsifier $\size S$ is linear in $k$ (with no logarithmic losses). It is a very challenging open problem to improve the bound of \cite{KPS24} to $\size S \leq O(k/\eps^2).$
\item 
Even though the graph sparsification problem has been extended to far more complex settings such as the hypergraph sparsification \cite{KKTY21,Lee23} or submodular functions sparsification \cite{JLLS24,Quan24}, the only known constructions of sparsifiers with no logarithmic loss is \cite{BSS14,MSS15}.
That being said, it is even more challenging to construct unweighted sparsifiers as essentially the only known construction is the solution to the notorious Kadison-Singer problem which was recently resolved in the breakthrough work of Marcus-Spielman-Srivastava \cite{MSS15}. As we see below the $L_1$ version of unweighted sparsifiers with no logarithmic dependency corresponds to long-standing open problems in graph theory namely the thin tree problem.
\item 
To the best of our knowledge, prior to our work the existence of unweighted code sparsifiers  was not studied (or perhaps expected).
\end{itemize} 

\paragraph{Applications to Graph Theory.}
We also briefly explain applications of this theorem to thin subgraphs. 
Given an unweighted (undirected) graph $G=(V,E)$ we say a set $T\subseteq E$ is a $\alpha$-thin w.r.t. $G$ if for any nonempty set $S\subsetneq V$,
$$ \abs{T(S,\overline{S})} \leq \alpha \abs{E(S,\overline{S})},$$
i.e. $T$ has at most $\alpha$-fraction of the edges of every cut. Recall that a graph $G$ is $\gk$-edge-connected if every cut in $G$ has at least $\gk$ edges.
The following thin tree conjecture is proposed by Goddyn two decades ago \cite{God04} and has been a subject of intense study since then \cite{AGMOS10,OS11,HO14,AO15,MP19,Mou19,Alg23, KO23}.
\begin{conj}[Thin Tree Conjecture]\label{conj:thintree}
    For any $\alpha<1$, there exists $\gk\geq 1$ such that any $\gk$-edge-connected graph $G$ has a spanning tree $T$ that is $\alpha$-thin. 
\end{conj}
We remark that there has been a "spectral" constructions of linear-sized thin subsets, see, e.g., \cite{BSS14,Ove15} but it remained an open problem whether one can construct thin subsets combinatorially without appealing to linear algebraic arguments. This question is specially motivated to address the thin tree conjecture since $\gk$-edge-connected graphs do not necessarily have spectrally thin trees, see \cite{AO15}. So, a combinatorial proof of the existence of thin subsets may potentially lead to a resolution of the thin tree conjecture. 

The following corollary is an application of \cref{thm:main}, \cref{cor:smallS} and \cref{cor:bigAlpha} to the cut-space of a graph.

\begin{cor}
    Let $G = \pair V E$ be a graph. Then 
    \begin{enumerate}
        \item $G$ has at least $2^{\size E - (\size V - 1)}$ many $1/2$-thin subgraphs. In other words, a uniformly random set $S\subseteq E$ is $1/2$-thin with probability at least $2^{1-\size V}$. \label{item:bigProb}
        \item $G$ has a $1/2$-thin subgraph $S$ with $\size S \geq \frac{\size E}{2} - \sqrt{\frac{\ln 2}{2} \size E \cdot (\size V - 1)}$. More generally, for any positive integer $\ell$, $G$ has a $2^{-\ell}$-thin subgraph with at least $\frac{\size E}{2^{\ell}} - c\sqrt{\size E(\size V-1)}$ edges for $c = \sqrt{\ln2}(\sqrt{2}+1)$.
        \label{item:smallThin}
    \end{enumerate}
\end{cor}

\begin{proof}
    Let $\cC \subseteq \bF_2^E$ denote the cut-space of the connected graph $G$. That is, $\cC$ contains all the vectors in $\bF_2^E$ that correspond to cuts in $G$. It is well-known that $\cC$ is a subspace of $\bF_2^E$ and that its dimension is $\size V - r$ where $r$ is the number of connected components of $G$. Thus in particular the dimension of $\cC$ is at most $\size V - 1$.
    \\
    Note that the complement of an $\alpha$-thin subgraph of $G$ corresponds to an unweighted one-sided $(1 - \alpha)$-sparsifier of $\cC$. Therefore \ref{item:bigProb} follows from \cref{thm:main} while  \ref{item:smallThin} follows from \cref{cor:smallS} and \cref{cor:bigAlpha}.
\end{proof}

We remark that although the existence of linear sized $1/2$-thin subgraphs (part \ref{item:smallThin} of the above) was known (up to worst constants) by spectral arguments such as \cite{BSS14,Ove15}, we are not aware of any inverse exponential in $\size V$ (independent of $\size E$) lower-bound on the probability that a uniformly random subgraph is $1/2$-thin.

\section{Main Proof}
Here we prove \cref{thm:main}. As in the statement of the theorem, let $\cC \subseteq \bF_2^n $ be a $\cd$-dimensional linear code.
\\
In what follows let us identify a set $S \subseteq [n]$ with its indicator vector $\mathbf{1}_S \in \bF_2^n$. Thus we identify $2^{[n]}$ with $\bF_2^n$.

\begin{defn}[codeword flip]
Given a codeword $c\in \C$, the flip corresponding to $c$ is the map $\bF_2^n\to \bF_2^n$ sending $S\in \mathbb{F}_2^n$ to $S+c$. 

When $\C$ is linear, this defines an equivalence relation on subsets of $[n]$. We say $S \sim S'$ if $S$ can be obtained from $S'$ by a codeword flip, i.e., 
$S + S' \in \C$. The transitivity of this relation follows from the linearity of $\C$. 
\end{defn}
Note that the collection of equivalence classes is the quotient space $\bF_2^n/\C$ (this is a vector space over $\bF_2$). Since $\dim \bF_2^n/\C = n - \cd$ there are precisely $\abs{\bF_2^n/\C}=2^{n-\cd}$ equivalence classes. 

\begin{lemma}
    \label{lemma:thinChar}
    Fix an equivalence class $H \in \bF_2^n/\C$. Let $S^*=\underset{S\in H}{\amax} \abs{S}$ 
    be a set with the largest size in $H$. 
    Then $S$ is an unweighted one-sided $1/2$-sparsifier for $\cC$.
\end{lemma}
\begin{proof}
    We prove this by contradiction. Suppose that there exists a codeword $c\in \C$ such that 
    $$ \wt(c_{S^*}) < \wt(c)/2.$$
    Let $S =c+S^*$. By definition $S\in H$. But the above equation implies $\abs{S}>\abs{S^*}$ which is a contradiction.   
\end{proof}
We are finally ready to prove \cref{thm:main}.
\begin{proof}[proof of \cref{thm:main}]
For an equivalence class $H$, let $S_H \subset [n]$ be a set of the largest size in $H$.
\Cref{lemma:thinChar} implies that the collection $\set{S_H \colon H \in \bF_2^n/\C}$ consists of $2^{n-\cd}$ distinct unweighted one-sided $1/2$-sparsifiers of $\cC$.
\end{proof}

\paragraph{Discussion.}
We say that a set $S\subseteq [n]$ is a hitting-set of a code $\C\subseteq \bF_2^n$ if for every (nonzero) codeword $c\in \C$, $c_S\neq 0$ i.e. every codeword has at least one coordinate in $S$. 
We conclude by proposing the following conjecture that generalizes the thin tree conjecture \ref{conj:thintree} to binary codes.
\begin{conj}
    For any $\alpha < 1$, there exists an integer $d = d(\alpha) \geq 1$ such that for every linear code $\cC \subseteq \bF_2^n$ with at least $d$ disjoint hitting sets, $S_1,\dots,S_d$, there exists a set $S\subseteq [n]$ such that for any (nonzero) $c\in \C$, $\alpha\cdot \wt(c)\leq \wt(c_S)<\wt(c)$.
\end{conj}
The above conjecture generalizes \autoref{conj:thintree} since every $d$-edge-connected graph has at least $d/2$ edge disjoint spanning trees by the Nash-Williams theorem.
We remark that a stronger version of this conjecture  (in an earlier version of this manuscript) was refuted by Putterman when applied to random linear codes.  

\appendix
\section*{Appendix}

\begin{proof}[proof of \cref{cor:bigAlpha}]
    The idea is as follows. Let $\gamma = \sqrt{\frac{\ln 2}{2}}$. 
    Let $S_1$ be the unweighted one-sided $1/2$-sparsifier of $\cC$ given by  \cref{cor:smallS} i.e. with $\size {S_1} \leq n/2 + \gamma \sqrt{nk}$. Without loss of generality, we may assume $\size{S_1} \geq n/2$ (as adding more coordinates does not hurt the sparsification property of $S_1$).
    \\
    Let $\cC_1'$ be the code given by projecting $\cC$ onto the coordinates in $[n]\setminus S_1$. Observe that $\cC_1'$ is a multiset as it might contain some codewords multiple times; let $\cC_1$ be the set of codewords that appears in $\cC_1'$ and observe that $\cC_1$ is a linear code of dimension at most $k$. Let $n_1 = n - \size {S_1}$ be the length of the new code.
    \\
    Now repeat this process with the new code $\cC_1$ for $\ell-1$ times. This process yields a sequence of sets $S_1,\ldots,S_\ell \subseteq [n]$ and codes $\cC_1,\ldots,\cC_\ell$ of lengths $n_1,\ldots,n_\ell$ satisfying the following for all $1 \leq r \leq \ell$ (let also define $n_0 = n, ~ \cC_0 = \cC$ for convenience).
    \begin{enumerate}
        \item $C_r$ is defined on coordinates in $[n] \setminus \U_{1 \leq i \leq r} S_r$.
        \item $n_r = n - \sum_{1 \leq i \leq r} \size {S_i}$. \label{item:nrec}
        \item $S_r \subseteq [n] \setminus \U_{i \leq r-1} S_i$.
        \item $\frac12 n_{r-1} \leq \size {S_r} \leq \frac12 n_{r-1} + \gamma \sqrt{n_{r-1}\cdot k}$. \label{item:Sub}
        \item $S_r$ is an unweighted one-sided $1/2$-sparsifier for $\cC_{r-1}$. \label{item:sparsifier}
    \end{enumerate}
    Note that the existence of a sparsifier $S_r$ satisfying the RHS of the inequality \ref{item:Sub} is guaranteed by \cref{cor:smallS}, the LHS inequality comes for free as adding more coordinates into a sparsifier does not hurt the sparsification property.
    \\
    For $0 \leq r \leq \ell$ let $T_r = [n] \setminus \U_{1 \leq i \leq r} S_i$.
    Note that $\cC_r$ is defined over coordinates in $T_r$ and hence by property \ref{item:sparsifier} we have $\wt(c_{T_r}) \leq \frac{1}{2}\wt(c_{T_{r-1}})$ thus $\wt(c_{T_r}) \leq \frac{1}{2^r}\wt(c)$ for all $c \in \cC$.
    \\
    Therefore letting $S = [n] \setminus T_\ell = \U_{1 \leq r \leq \ell} S_r$, we have that $S$ is an unweighted one-sided $(1- \frac{1}{2^\ell})$-sparsifier for $\cC$.
    \\
    We now provide the upper bound for $\size S = n - n_\ell$. 
    Properties \ref{item:nrec} and \ref{item:Sub} imply 
    \[
        \frac{1}{2}n_{r-1} - \gamma \sqrt{n_{r-1} \cdot k}
        \leq
        n_r 
        \leq
        \frac{1}{2}n_{r-1}
        .
    \]
    Note that since $n_r \leq \frac 12 n_{r-1}$ it follows inductively that $n_r \leq \frac{1}{2^r}n$ for all $0 \leq r \leq \ell$. Therefore
    \[
        \frac{1}{2}n_{r-1} - \gamma \sqrt{\frac{n}{2^{r-1}}k}
        \leq
        n_r.
    \]
    Hence inductively
    \[
        n_r \geq \frac{1}{2^r}n - \gamma \sqrt{nk} \cdot \sum_{i = 0}^{r-1} 2^{-i/2} \geq \frac{n}{2^r} - c\sqrt{nk},
    \]
    for $c=\gamma\frac{1}{1-1/\sqrt{2}}$ as defined in the statement of the corollary.
    This implies $\size S = n - n_\ell\leq (1 - \frac{1}{2^\ell})n + c \sqrt{nk} = \alpha n + c \sqrt{nk}$, as desired.
\end{proof}

\printbibliography

\end{document}